\documentclass[11pt]{article}
\usepackage[dvips]{graphicx}
\usepackage{epsfig}
\usepackage{multirow}
\usepackage{amsthm}
\usepackage{amscd}
\usepackage{amsmath}
\usepackage{amsfonts}
\usepackage{amssymb}
\usepackage{multirow}
\usepackage[authoryear, comma]{natbib}
\usepackage{algorithm}
\usepackage{algorithmic}
\newtheorem{theo}{Theorem}
\newtheorem{prop}{Proposition}

\newtheorem{defin}{Definition}

\newtheorem*{rem}{Remark}

\def\one{{\bf 1}\hskip-.5mm}

\def\E{{\mathbb E}}
\def\F{{\mathcal F}}
\def\P{{\mathbb P}}

\def\Z{{\mathbb Z}}

 \def\F {{\mathcal F}}

	\makeatletter
	\renewcommand\@biblabel[1]{#1.}
	\makeatother

\begin{document}

\title{Identifying interacting pairs of sites in Ising models on a countable set}
\date{December 16, 2014}
\author{Antonio Galves, Enza Orlandi and Daniel Yasumasa Takahashi}
\maketitle

 \begin{center} {\sl Para Errico, com um grande abra\c{c}o} \end{center}

\begin{abstract}
This paper address the problem of identifying pairs of interacting  sites from a finite sample of independent realizations of the Ising model. We consider Ising models  in a infinite countable set of sites under Dobrushin uniqueness condition. The observed sample contains only the values assigned by the Ising model to a finite set of sites. Our main result is an upperbound for the probability of misidentification of the pairs of interacting sites in this finite set.



\end{abstract}

\section{Introduction}
%
%The class of Ising models is a popular family of exponential distributions on the set of binary configurations supported by a finite or a countable set of sites. The law of the Ising model is characterized by a pairwise potential, \textit{i.e.}, a family of real numbers $J(i,j)$ indexed by pairs of sites $(i,j)$. Originally introduced in statistical mechanics as a mathematical model for ferromagnetism, the Ising model has been extensively used, for instance, in computer vision \citep{Woods78, Besag93}, image processing \citep{Cross83}, neuroscience \citep{Schneidman06}, and as a general model in spatial statistics \citep{Ripley81}. The references given above are just starting points of a huge literature.   For a recent statistical physics oriented survey of rigorous mathematical results on  Gibbs distributions, including Ising models, we refer the reader to \citet{Presutti09}.

 In this article we address the statistical problem of identifying the pairs of interacting sites of an Ising model on a countable set of sites (possibly infinite) when we only observe the values assigned  on a  finite subset of sites (partially observed). Our sample consists of a finite number of independent realizations of the Ising model observed at this finite set. We introduce a statistical procedure to identify the interacting pair of sites given the observations. Our main result is an upper bound for the probability of misidentifying the pairs of interacting sites. 
% 
% The law of the Ising model is characterized by a pairwise potential, \textit{i.e.}, a family of real numbers indexed by pairs of sites representing the interaction weights between pair of sites. We are interested on identifying the pair of sites with $(i,j)$ for which $J(i,j) \neq 0$. 
 
 Originally introduced in statistical mechanics as a mathematical model for ferromagnetism, the Ising model has been extensively used, for instance, in computer vision \citep{Woods78, Besag93}, image processing \citep{Cross83}, neuroscience \citep{Schneidman06}, and as a general model in spatial statistics \citep{Ripley81}. The references given above are just starting points of a huge literature.   For a recent statistical physics oriented survey of rigorous mathematical results on  Gibbs distributions, including Ising models, we refer the reader to \citet{Presutti09}.

When the set of sites is not finite, the Ising model is supported by a continuous set of infinite configurations. However, from an applied statistics point of view, we cannot observe more than the projection of the Ising model on a finite subset of sites. We introduce an estimator for the set of interacting pairs of sites belonging to the finite set we observe.
This estimator can be informally described as follows. For each site $i$ in the observed finite set we estimate the conditional probability of the model in $i$, given the remaining sites in the finite set. Then we compare this empirical conditional probability with the empirical conditional probability on the same site $i$ given the remaining sites with the exception of another site $j$, with $j\neq i$. If the two conditional probabilities are statistically equal, we conclude that interacting weight between the sites is null. 

 The proof of the main theorem has two ingredients, which are interesting by themselves. 
The first ingredient is an upperbound for the probability of misidentification for the Ising model on a finite set of sites. This is the content of Theorem \ref{teo:bounded}. 
The second ingredient is a coupling result given in Theorem \ref{teo:dob}. It says that we can couple together an Ising model on infinite set of sites and an Ising model restricted to a finite subset of sites in such a way that the probability of a discrepancy at a fixed site vanishes as the set of observed sites increases to the entire set. As a consequence of this result, we are able to bound above the probability of misidentification due to the fact that we are able to observe only a finite set of sites, not the entire set of interacting sites.  The proof of this result uses a constructive version of classical Dobrushin's contraction method. For a nice presentation of the contraction method in its original framework we refer the reader to \citet{Presutti09}.

It is important to note that we don't need to assign a metric to the set of sites to state and prove our results. It turns out that in several applied contexts, a predefined metric is unwarranted.  For instance, for the problem of inferring the presence/absence of interactions between pairs of neurons, it is not clear \emph{a priori} that there is any  consistent relationship between the strength of the interaction and the physical distance (or any other metric) between the neurons. Moreover, in several situations, the experimenter doesn't know if the recorded neurons are physically close or far.  This justifies the approach taken here.

Let conclude this section with short comments on the recent related literature. The case of random field on a finite set of sites, which is entirely observed  was considered in \citet{Ravikumar10, Bento09, Bresler08}. The infinite case was also considered in \cite{locherbach2011neighborhood, Csiszar06, lerasle2014sharp, lerasle2011oracle}. The first two use a BIC like approach for an homogeneous random field using a single observation and the last two use an oracle approach, which solves a problem that is different from the identification problem considered in this article.
 
This paper is organized as follows.   Notation, definitions and main results are presented in Section \ref{sec:results}. The proofs of the theorems are presented in Section \ref{sec:proof}.

\section{Notation, definitions, and main results} \label{sec:results}

Let $S$ be a countable set of sites. 
\begin{defin}
  A pairwise potential is a family $J=\{ J(i,j): (i, j) \in S
  \times S \} $ of real numbers which satisfy the conditions
\begin{equation}
J(i,i) = 0, \qquad   J(i,j)=J(j,i), \qquad  \sup_{i \in S} \sum_{j \in S}   |J(i,j)|  <\infty .  
\end{equation} 
\end{defin}

  Let $\mathcal{X}=\{-1,1\}^{S} $ be the set of configurations on the set of sites $S$.  A fixed configurations will be denoted by lower case letter $x$ whereas the capital letter $X$ will denote a random configuration taking values on $\mathcal{X}$ and probability measure $\P$.   For any $i \in S$, $x(i)$ will denote the value of the configuration $x$ at site $i$. Given a subset $F$ of $S$, we shall also denote $x(F)=\{x(i): i \in F\}$ and similarly for $X$. 
%We also denote by $V_i$ the set of sites $j \in V$ such that $(i,j) \in E$.

\begin{defin}
The Ising model with pairwise potential $J$ is a random configuration $X$ with values on $\mathcal{X}$, which probability satisfies 
\begin{equation*}
  \mathbb{P} \left( X (i) = x (i) | X (j)=x (j)\, ,\, j\neq i \right) = 
\frac{1}{1 + \exp ( -2 \sum_{j\in S} J(i,j) x (i) x (j) )},
\end{equation*}
for all $i \in S$ and for $\P$-a.e. $x \in \mathcal{X}$. 
\end{defin}
In the above definition the left hand side of the above equality denotes a regular version of the conditional
probability of $X (i)$ given that $X(j)=x(j)$ for $ j\neq i$. 

Let $F$ be a finite subset of $S$.  We  use the shorthand notation $p(x(F))$ and 
$p( x(i) | x(F) )$ to denote, respectively, the probability $\mathbb{P}(X(F)=x(F))$ and the conditional probability $\mathbb{P}(X(i)=x(i) | X(F)=x(F ))$.

\begin{defin}
For any site $i\in S$, the interaction neighborhood $\mathcal{G}(i)$ is defined as
\begin{equation*}  
\mathcal{G}(i)=\left\{ j \in S : J(i,j) \neq 0 \right\}.
\end{equation*}
\end{defin}

In general, we cannot observe the entire random configuration on $S$, but only the values on some finite subset $F \subset S$. Moreover, we observe only a finite number of samples, \emph{i.e.}, the observations are i.i.d. samples $X_1(F),\ldots,X_n(F)$. In this situation, we do not expect to recover $\mathcal{G}(i)$, but we might be able to identify $\mathcal{G}(i) \cap F$. In this article, we show how we can do it. The following family of sets will be useful for the rest of the article.

\begin{defin}
 A family $\mathcal{F}$ of finite subsets $F_i \subset S$ indexed by $i \in S$  is called a truncation class if for any $i,j \in S$ we have that $i \in F_i$ and  $j \in F_i \iff i \in F_j$.
\end{defin}

It is convenient to introduce the following truncated version of the Ising model.

\begin{defin}
Given a truncation class $\mathcal{F}$, we denote by $J^\mathcal{F}$ the truncated potential defined as follows
\begin{equation} \label{trunJ}
J^\mathcal{F}(i,j)=
\begin{cases} 
 J(i,j)\; & \mbox{if } j \in F_i \\
  0\;, & \mbox{otherwise} ,
\end{cases}
\, 
\end{equation}
We also denote by $X^\F$ the corresponding Ising model with pairwise potential $J^\F$. 
\end{defin}

From now on $\F$ will always denote a truncation class. Let $F$ be a finite subset of $S$. As before, we use the shorthand notation $p^\F(x(F))$ and 
$p^\F( x(i) | x(F) )$ to denote, respectively, the probability $\mathbb{P}(X^\F(F)=x(F))$ and the conditional probability $\mathbb{P}(X^\F(i)=x(i) | X^\F(F)=x(F))$.

Given a site $i \in S$,  and a finite set $F_i \in \F$, let
$$D(x,F_i,i,j) = \left |p^\F( x(i)|x(F_i\setminus \{i\}))-p^\F\left(x(i)|x(F_i\setminus \{i,j\})\right)\right|p^\F(x(F_i\setminus \{i\})). $$

\begin{defin} \label{def:rerldif}

 For any $i \in S$, $F_i \in \F$, and $\epsilon > 0$, the interaction neighborhood $V^\F(i)$ of $i$ is defined as 

\begin{equation}  \label{eq:rerldif}
V^\F(i)=\left\{ j \in F_i : \max_{x(F_i)} D(x,F_i,i,j)  > 2\epsilon \right\}.
\end{equation}
\end{defin}

 We prove the following lowerbound.
\begin{prop}  \label{prop:lowerbound}
 Let $i \in S$, $F_i \subset \F$, $j \in F_i$, and $\sup_{i \in S} \sum_{j \in S}   |J(i,j)| = \gamma$. We have that
\begin{equation}
 \max_{x(F_i)} D(x,F_i,i,j) \geq \frac{2^{-|F_i|+2}e^{2\gamma}}{(1+e^{2\gamma})}|J(i,j)| 
 \end{equation}
\end{prop}
Observe that if $j \in \mathcal{G}(i) \cap F_i$, we have $|J(i,j)| > 0$, therefore for small enough $\epsilon$, we will have $V^\F(i) = \mathcal{G}(i) \cap F_i$.

Let us introduce, for any finite set $F$, the empirical probability measure
\begin{equation*}
\hat{p}_n(x(F) )=\frac{1}{n} \sum_{k=1}^n\one\{X_k(F)=x(F)\},
\end{equation*}
where $\one$ denotes the indicator function. Given any site $j \in S$, we will also define the empirical conditional probability
$$\hat{p}_n(x(j) | x(F\setminus \{j\} ))=\frac{\hat{p}_n(x(F \cup \{j\}) )}{\hat{p}_n(x(F \setminus \{j\}) )}\; ,$$ 
if $\hat{p}_n(x(F \setminus \{j\}) ) > 0 $ and $ \hat{p}_n(x(j) | x(F\setminus \{j\} )) = 0$, otherwise.

For any $i \in S$, $F_i \in \F$, any configuration $x(F_i)$, and $j \in F_i$ we define the empirical weighted distance between the conditional probabilities as follows
\begin{equation}\label{eq:difference}
\hat{D}_n(x,F_i,i,j)= 
\left|\hat{p}_n(x(i) | x(F_i\setminus\{i\}))-\hat{p}_n(x(i) | x(F_i\setminus\{i, j\}))\right|
\hat{p}_n(x(F_i\setminus \{i\})\,.
\end{equation}
Note that $\hat{D}_n(x,F_i,i,j)$ is a function of the sample $X_1,\ldots,X_n$ and is therefore a random variable.

We can now define our estimator.
\begin{defin} \label{def:estimatornei}
 For any $i \in S$ and $F_i \in \F$, the interaction neighborhood estimator  is  defined as
\[
\hat{V}_n(i)=\left\{ j \in F_i : \max_{x(F_i)} \hat{D}_n(x,F_i,i,j)  >  \epsilon \right\}\, ,
\]
where the threshold $\epsilon$ is the same as in Definition \ref{def:rerldif}.
\end{defin}

We can now state our main result. 

\begin{theo} \label{teo:main}
Let  $i \in S$, $F_i \in \F$, and $X_1(F_i), \ldots, X_n(F_i)$ be the local projections of independent realizations of an Ising model whose pairwise potential satisfies 
\begin{equation} \label{cond:J}
\sup_{k\in S}\sum_{j\in S} |J(k,j)|=r<1.
\end{equation} 
Then, for any threshold value $\epsilon > 0$, we have
\begin{align}
&\mathbb{P}\left(\hat{V}_n(i) \neq V^\F(i) \right) \notag \\
& \le 4\exp\left(-\frac{n\epsilon^2}{8v+\frac{4}{3}\epsilon}+2|F_i| \right) + \frac{1}{1-r}n|F_i|\left(\sup_{k \in S}\sum_{ j \in S \setminus F_k} |J(k,j)|\right)\,,\label{eq:uppermain1}
\end{align}
where 
\begin{equation}
v=\sup_{x(F_i)}\sup_{j \in F_i}\left(1-p^\F\left(x(i)|x(F_i\setminus \{i,j\})\right)p^\F(x(F_i))\right).  \label{eq:variancemain}
\end{equation}
\end{theo}

\vspace{0.5cm}

The condition (\ref{cond:J}) is known as Dobrushin uniqueness condition in the statistical physics literature \citep{Presutti09}.

\vspace{1cm}

The first ingredient in the proof of  Theorem \ref{teo:main} is an upperbound for the probability of misidentification of interacting pairs in the case of a finite range interaction.   This is given in the next theorem.

\begin{theo} \label{teo:bounded}
Let $i \in S$,  $F_i \in \F$,  and $X_1^\F(F_i), \ldots, X_n^\F(F_i)$ be the projections of independent realizations of an Ising model with pairwise potential $J^\F$.

Then for any site $i \in S$ and any threshold value $\epsilon > 0$, we have
\[
\mathbb{P}\left(\hat{V}_n(i) \neq V^\F(i) \right) 
\le 4 \exp\left(-\frac{n\epsilon^2}{8v+\frac{4}{3}\epsilon} + 2|F_i| \right)\,,
\] 
where 
\begin{equation*}
v=\sup_{x(F_i)}\sup_{j \in F_i}(1-p^\F\left(x(i) | x(F_i\setminus \{i,j\})\right)p^\F(x(F_i))). 
\end{equation*}
\end{theo}

The second ingredient in the proof of Theorem \ref{teo:main} is a coupling result. To state it we first need to introduce the definition of coupling.

\begin{defin}
Let $X$ and $X^\F$ be Ising models with pairwise potentials $J$ and $J^\F$, respectively. A coupling between $X$ and $X^\F$ is a random element $(\tilde{X} , \tilde{X}^\F )$ taking values on $S\times S$ such that
\begin{enumerate}
\item $\tilde{X}$ has the same law as $X$;
\item $\tilde{X}^\F$ has the same law as $X^\F$.
\end{enumerate}
\end{defin}

The following theorem says that we can sample together $X$ and $X^\F$ and gives an upper bound for the probability of discrepancy between $X(i)$ and $X^\F(i)$.

\begin{theo} \label{teo:dob}
If $J$ is pairwise potential which satisfies condition (\ref{cond:J}), \textit{i.e.},
\begin{equation} \label{eq:dob}
\sup_{k \in S}\sum_{j \in S}|J(k,j)| = r < 1 
\end{equation}
and $J^\F$ is defined as in (\ref{trunJ}), then there exists a coupling $(\tilde{X} , \tilde{X}^\F )$ such that for any $i \in S$ the following inequality holds
\begin{equation} 
\P\left(\tilde{X}(i) \neq \tilde{X}^\F(i) \right) \leq \frac{1}{1-r} \sup_{k\in S}\sum_{j \in S \setminus F_k}|J(k,j)|. \label{bound2}
\end{equation}

\end{theo}

\vspace{0.5cm}
\begin{rem} In practice, it is important to find a truncation class $\F$  that makes the right hand side of  (\ref{bound2}) small. A simple example is given by an Ising model $X$ with $S = \Z^d$ and nearest neighborhood interaction, \textit{i.e.}, $\mathcal{G}(i) = \{j \in \Z^d: |i-j| = 1\}$. The right hand side of (\ref{bound2}) will be zero if we take $\F$ such that, for all $k\in S$, $\mathcal{G}(k) \subset F_k$.
\end{rem}

\vspace{1cm}
\section{Proof of the results} \label{sec:proof}

\section*{Proof of Proposition \ref{prop:lowerbound}}

Let $i, j\in F_i$, $x(F_i\setminus \{i\}) \in\{-1,+1\}^{|F_i|-1} $, and  $y(F_i\setminus \{i\}) \in\{-1,+1\}^{|F_i|-1} $ with $y(j) = -x(j)$. Using the mean value theorem, we have
\begin{equation*}
\left|p^\F( x(i)|x(F_i\setminus \{i\}))-p^\F\left(x(i)|y(F_i\setminus \{i\})\right) \right|
 \geq \frac{e^{2\gamma}}{(1+e^{2\gamma})^2}|J(i,j)|.
\end{equation*}

Hence, for any $j \in F_i$ such that $J(i,j) \neq 0$ we have
\begin{align*}
&\left|p^\F( x(i)|x(F_i\setminus \{i\}))-p^\F\left(x(i)|x(F_i\setminus \{i,j\})\right) \right| \\
& \geq \frac{2e^{2\gamma}}{(1+e^{2\gamma})^2}|J(i,j)|\min_{x(j) \in \{-1,+1\}}  p^\F(x(j) | x(F_i\setminus \{i,j\})).
\end{align*}
Also
$$\min_{x(j) \in \{-1,+1\}}  p^\F(x(j) | x(F_i\setminus \{i,j\})) \geq \frac{1}{1+e^{2\gamma}}.$$

We observe that for any $i \in F$ 
\begin{equation*}
\max_{x(F_i\setminus \{i\})} p^\F(x(F_i\setminus \{i\})) \geq 2^{-|F_i|+1}.
\end{equation*}

Combining the above inequalities, we have
\begin{equation*}
\max_{x(F_i)} D(x,F_i,i,j) \geq \frac{2^{-|F_i|+2}e^{2\gamma}}{(1+e^{2\gamma})^3}|J(i,j)|,
\end{equation*}
as we wanted to show.

\qed

\section*{Proof of Theorem \ref{teo:main}}

Let the finite set $F \subset S$ be fixed and let $(\tilde{X}_1,\tilde{X}_1^\F), \ldots, (\tilde{X}_n,\tilde{X}_n^\F)$ be $n$ independent copies of the pair $(\tilde{X}, \tilde{X}^\F)$ which existence is guaranteed by Theorem \ref{teo:dob}.  The random elements $\tilde{X}_1, \ldots, \tilde{X}_n$ are independent copies of the Ising model $X$ with pairwise potential $J$.  The random elements $\tilde{X}^\F_1, \ldots, \tilde{X}^\F_n$ are independent copies of the Ising model $X^\F$ with truncated pairwise potential $J^\F$ defined as in (\ref{trunJ}).  

Let us indicate explicitly the sample in all the statistics and events appearing in Theorem  \ref{teo:main} as functions either of the sample $\tilde{X}_1, \ldots, \tilde{X}_n$ or of the sample $\tilde{X}^\F_1, \ldots, \tilde{X}^\F_n$ . We start with notation of the empirical probability measures $\hat{p}_n$, as follows
\begin{align*}
&\hat{p}_n(x(F_i) )[\tilde{X}_1, \ldots, \tilde{X}_n]=\frac{1}{n} \sum_{k=1}^n\one\{\tilde{X}_k(F_i)=x(F_i)\}\\
&\hat{p}_n(x(F_i) )[\tilde{X}^\F_1, \ldots, \tilde{X}^\F_n]=\frac{1}{n} \sum_{k=1}^n\one\{\tilde{X}^\F_k(F_i)=x(F_i)\}.
\end{align*}

To simplify the writing we shall use the short notation $$\tilde{\mathbf{X}}=(\tilde{X}_1, \ldots, \tilde{X}_n)\;\;\;\text{and}\;\; \;\tilde{\mathbf{X}}^\F=(\tilde{X}^\F_1, \ldots, \tilde{X}^\F_n).$$

Now using either the empirical probability measures $\hat{p}_n(x(F_i) )[\tilde{\mathbf{X}}]$ or $\hat{p}_n(x(F_i) )[\tilde{\mathbf{X}}^\F]$ we define the neighborhood estimators $\hat{V}_n(i)[\tilde{\mathbf{X}}]$ and $\hat{V}^\F_n(i)[\tilde{\mathbf{X}}^\F]$.

Now we are ready to conclude the proof.
An upperbound for the probability of misidentification for the sample $\tilde{X}_1, \ldots, \tilde{X}_n$ is given by  
\begin{align*}
&\P\left(\hat{V}_n(i)[\tilde{\mathbf{X}}] \neq V^\F(i) \right)\\
&\leq \P\left(\left\{\hat{V}_n(i)[\tilde{\mathbf{X}}] \neq V^\F(i) \right\} \bigcap \bigcap_{k\in \{1,\ldots, n\}} \bigcap_{j\in F_i} \left\{\tilde{X}_{k}^{\F}(j)=\tilde{X}_{k}(j)\right\}\right)\\
&+\P\left(\bigcup_{ k\in \{1,\ldots, n\}} \bigcup_{j\in F_i}\left\{\tilde{X}_{k}^{\F}(j) \neq \tilde{X}_{k}(j) \right\} \right).
\end{align*}

By Theorem \ref{teo:dob}
\begin{equation*}
\P\left( \bigcup_{ k\in \{1,\ldots, n\}} \bigcup_{j\in F_i}\left\{\tilde{X}_{k}^{\F}(j) \neq \tilde{X}_{k}(j) \right\} \right)\\
\leq n|F_i|\sup_{k \in S} \sum_{ j \in S\setminus F_k} |J(k,j)|.
\end{equation*}

Now, we observe that in the set 
$$ \bigcap_{k\in \{1,\ldots, n\}} \bigcap_{j\in F} \left\{\tilde{X}_{k}^{\F}(j)=\tilde{X}_{k}(j)\right\}$$
the following holds
$$\hat{V}_n(i)[\tilde{\mathbf{X}}] =\hat{V}_n(i)[\tilde{\mathbf{X}}^\F].$$

Hence
\begin{align}
 &\P\left(\left\{\hat{V}_n(i)[\tilde{\mathbf{X}}] \neq V^\F(i) \right\}  \bigcap \bigcap_{k\in \{1,\ldots, n\}} \left\{\tilde{X}_{k}^{\F}(F_i)=\tilde{X}_{k}(F_i)\right\}\right) \notag\\
 &=\P\left(\left\{\hat{V}_n(i)[\tilde{\mathbf{X}}^\F] \neq V^\F(i) \right\}  \bigcap \bigcap_{k\in \{1,\ldots, n\}} \left\{\tilde{X}_{k}^{\F}(F_i)=\tilde{X}_{k}(F_i)\right\}\right) \notag\\
 &\leq \P\left(\hat{V}_n(i)[\tilde{\mathbf{X}}^\F] \neq V^\F(i)  \right) . \label{eq:probultimo}
 \end{align}
 
 Since Theorem \ref{teo:bounded} provides an upperbound for the last term in (\ref{eq:probultimo}), we have
 \begin{align*}
 &\P\left(\left\{\hat{V}_n(i)[\tilde{\mathbf{X}}] \neq V^\F(i) \right\}  \bigcap \bigcap_{k\in \{1,\ldots, n\}} \left\{\tilde{X}_{k}^{\F}(F_i)=\tilde{X}_{k}(F_i)\right\}\right) \\
 &\leq 4\exp\left(-\frac{n\epsilon(F_i,n)^2}{8v+\frac{4}{3}\epsilon(F_i,n)}+2|F_i| \right).
\end{align*}
This concludes the proof of Theorem \ref{teo:main}.

\qed

% To prove Corollary \ref{cor:2} we observe that the finite volume case is equivalent to consider the pairwise potential satisfying
% \begin{equation} 
%J(i,j)=0\;\; \text{for} \;\; j \in \Z^d \setminus V
%\end{equation} 
%  in the $\Z^d$ case.

 \section*{Proof of Theorem \ref{teo:bounded} }

For convenience of the reader, before the proof let us recall the classical inequality of Bernstein which will be used in the sequence.

\noindent {\bf Bernstein inequality}
Let $\xi_1, \ldots, \xi_n$ be i.i.d. random variables with $|\xi_1| \leq b$ a.s. and $\E[\xi_1^2] \leq v < \infty$. Then the following inequality holds
\begin{equation*}
\P\left(\left|\frac{1}{n}\sum_{k=1}^n\xi_k - \E[\xi_1] \right| \geq \epsilon \right) \leq 2 \exp \left(-\frac{n \epsilon^2}{2(v+\frac{1}{3}b\epsilon)}\right).
\end{equation*}
For a proof of this inequality, we refer the reader to \citet{Massart03}.

To begin the proof of Theorem \ref{teo:bounded}, let us denote 
\begin{equation}\label{fpositive}
\mathcal{O}_n^\F(i)=\left\{ j \in \hat{V}_n^\F(i): j\in F_i\setminus V^\F(i) \right\}
\end{equation}
the event of {\sl false positive identification}. 

The event of {\sl false negative identification} is defined as
\begin{equation}\label{fnegative}
\mathcal{U}^\F_n(i)=\left\{j \in F_i \setminus \hat{V}_n(i): j\in  V^\F(i) \right\}.
\end{equation}

We observe that 
$$\{\hat{V}_n(i) \neq V^\F(i)\} = \mathcal{O}^\F_n(i) \cup \mathcal{U}^\F_n(i).$$

We will first  obtain an upperbound for the probability of event false positive identification.
Observe that 
\begin{equation} \label{eq:over1}
\mathbb{P}\left(\mathcal{O}^\F_n(i) \right) \leq  \sum_{x(F_i)}\;\; \sum_{j \in  F_i\setminus V^\F(i)}  \mathbb{P}\left( \hat{D}_n(x,F_i,i, j) > \epsilon \right)\, .
\end{equation}

Let us fix $j \in F_i\setminus V^\F(i)$ and $x(F_i)\in \{-1,+1\}^{F_i}$. 
To obtain an upperbound for the right hand side of (\ref{eq:over1}) we first observe that
\begin{align}
&\hat{D}_n(x,F_i,i,j) \notag\\
&\leq \left|\hat{p}_n(x(F_i))-p^\F\left(x(i)|x(F_i\setminus \{i,j\})\right)\hat{p}_n(x(F_i\setminus\{i\}))\right|\notag\\
&+\left|\frac{\hat{p}_n(x(F_i\setminus \{j\}))}{\hat{p}_n(x(F_i\setminus \{i,j\}))}-p^\F\left(x(i)|x(F_i\setminus \{i,j\})\right)\right|\hat{p}_n(x(F_i\setminus \{i\})). \label{eq:eqdup}
\end{align}
This inequality was obtained by adding and subtracting $$p^\F\left(x(i)|x(F_i\setminus \{i,j\})\right)\hat{p}_n(x(F_i\setminus\{i\}))$$ in expression (\ref{eq:difference}).

Since 
$$0\leq \frac{\hat{p}_n(x(F_i\setminus \{i\}))}{\hat{p}_n(x(F_i\setminus \{i,j\}))}\leq 1,$$
we finally obtain the upperbound
\begin{align}
&\hat{D}_n(x,F_i,i,j) \notag \\
&\leq \left|\hat{p}_n(x(F_i))-p^\F\left(x(i) | x(F_i\setminus \{i,j\})\right)\hat{p}_n(x(F_i\setminus \{i\}))\right|\notag\\
&+\left|\hat{p}_n(x(F_i\setminus \{j\}))-p^\F\left(x(i) | x(F_i\setminus \{i,j\})\right)\hat{p}_n(x(F_i\setminus \{i,j\}))\right|. \label{eq:upp1}
\end{align}

Therefore,
\begin{align*}
&\mathbb{P}\left( \hat{D}_n(x,F_i,i,j) > \epsilon \right)\\
& \leq \mathbb{P}\left( \left|\hat{p}_n(x(F_i))-p^\F\left(x(i)|x(F_i\setminus \{i,j\})\right)\hat{p}_n(x(F_i\setminus \{i\}))\right| > \frac{1}{2}\epsilon\right) \\
& +  \mathbb{P}\left( \left|\hat{p}_n(x(F_i\setminus \{j\}))-p^\F\left(x(i)|x(F_i\setminus \{i,j\})\right)\hat{p}_n(x(F_i\setminus \{i,j\}))\right|> \frac{1}{2}\epsilon\right).
\end{align*}

 The classical Bernstein inequality provides the following upperbounds for the terms in the right hand side of the above equation
\begin{align}
&\mathbb{P}\left( \left|\hat{p}_n(x(F_i))-p^\F\left(x(i)|x(F_i\setminus \{i,j\})\right)\hat{p}_n(x(F_i\setminus \{i\}))\right| > \frac{1}{2}\epsilon\right) \notag \\
&\leq 2\exp\left(-\frac{n\epsilon^2}{8v+\frac{4}{3}\epsilon} \right), \label{eq:firstbound}
\end{align}
where 
\begin{equation}
v=\sup_{x(F_i)}\sup_{j \in F_i}(1-p^\F\left(x(i)|x(F_i\setminus \{i,j\})\right)p^\F(x(F_i))). \label{eq:variance1}
\end{equation}
Also
\begin{align} 
&\mathbb{P}\left( \left|\hat{p}_n(x(F_i\setminus \{j\}))-p^\F\left(x(i)|x(F_i\setminus \{i,j\})\right)\hat{p}_n(x(F_i\setminus \{i,j\}))\right| > \frac{1}{2}\epsilon\right) \notag \\
&\leq 2\exp\left(-\frac{n\epsilon^2}{8v'+\frac{4}{3}\epsilon} \right), \label{eq:secbound}
\end{align}
where 
\begin{equation}
v'=\sup_{x(F_i)}\sup_{j \in F_i}(1-p^\F\left(x(i)|x(F_i\setminus \{i,j\})\right)p^\F(x(F_i\setminus \{j\}))). \label{eq:variance2}
\end{equation}

Summing up inequalities (\ref{eq:firstbound}) and (\ref{eq:secbound})  for all configurations $x(F_i)$ and all sites $j \in F_i\setminus V^\F(i)$ we obtain the following upperbound for the probability of false positive identification
\begin{align}
&\mathbb{P}\left(\mathcal{O}^\F_n(i) \right) \notag\\
&\leq 4 (|F_i|-|V^\F(i)|)\exp\left(-\frac{n\epsilon^2}{8v+\frac{4}{3}\epsilon} \right) \notag\\
&\leq 4 (|F_i|-|V^\F(i)|)\exp\left(-\frac{n\epsilon^2}{8v+\frac{4}{3}\epsilon}  \right). \label{eq:over}
\end{align}

We will now obtain an upperbound for the probability of false negative identification. 
For any $j \in V^\F(i)$ we have
\begin{equation} \label{eq:over}
\mathbb{P}\left(j \notin \hat{V}_n(i)\right) = \P\left( \bigcap_{x(F_i)} \left \{\hat{D}_n(x,F_i,i,j) \leq \epsilon  \right\}\right).
\end{equation}
To obtain an upperbound for (\ref{eq:over}), it is enough to obtain an upperbound for 
\begin{equation}
 \P\left(\hat{D}_n(x,F_i,i,j) \leq \epsilon \right) \label{eq:simplebound}
 \end{equation}
where $x(F_i)$ is any fixed configuration. In particular, we can take a configuration which maximizes
\begin{equation} \label{eq:threshold}
\left |p^\F( x(i)|x(F_i\setminus \{i\}))-p^\F\left(x(i)|x(F_i\setminus \{i,j\})\right)\right|p(x(F_i\setminus \{i\}).
\end{equation}

To do this, we first obtain a lower bound for $\hat{D}_n(x,F_i,i,j)$ in the same way we obtained the upperbound in (\ref{eq:upp1}).
\begin{align*}
&\hat{D}_n(x,F_i,i,j)\notag\\
&\geq \left|\hat{p}_n(x(F_i))-p^\F\left(x(i)|x(F_i\setminus \{i,j\})\right)\hat{p}_n(x(F_i\setminus\{i\}))\right|\notag\\
&-\left|\hat{p}_n(x(F_i\setminus \{j\}))-p^\F\left(x(i)|x(F_i\setminus \{i,j\})\right)\hat{p}_n(x(F_i\setminus \{i,j\}))\right|\frac{\hat{p}_n(x(F_i\setminus \{i\}))}{\hat{p}_n(x(F_i\setminus \{i,j\}))}.\notag
\end{align*}
Observing again that 
$$0\leq \frac{\hat{p}_n(x(F_i\setminus \{i\}))}{\hat{p}_n(x(F_i\setminus \{i,j\}))}\leq 1$$
we finally obtain the lower bound
\begin{align}
&\hat{D}_n(x,F_i,i,j) \notag\\
& \geq \left|\hat{p}_n(x(F_i))-p^\F\left(x(i)|x(F_i\setminus \{i,j\})\right)\hat{p}_n(x(F_i\setminus\{i\}))\right|\notag\\
&-\left|\hat{p}_n(x(F_i\setminus \{j\}))-p^\F\left(x(i)|x(F_i\setminus \{i,j\})\right)\hat{p}_n(x(F_i\setminus \{i,j\}))\right|. \label{eq:low2}
\end{align}

To make formulas shorter let us call for the moment
\begin{equation*}
 W = \hat{p}_n(x(F_i))-p^\F\left(x(i)|x(F_i\setminus \{i,j\})\right)\hat{p}_n(x(F_i\setminus\{i\}))
\end{equation*}
and
\begin{equation*}
 R = \hat{p}_n(x(F_i\setminus \{j\}))-p^\F\left(x(i)|x(F_i\setminus \{i,j\})\right)\hat{p}_n(x(F_i\setminus \{i,j\})).
\end{equation*}

With this new notation, using inequalities (\ref{eq:upp1}) and (\ref{eq:low2}) we obtain
\begin{equation}
\left| \hat{D}_n(x,F_i,i,j) - |W| \right| \leq |R |. \label{eq:enquadra}
\end{equation}

A straightforward computation shows that
\begin{equation*}
\mathbb{E}[W] = p(x(F_i))-p^\F\left(x(i)|x(F_i\setminus \{i,j\})\right)p^\F(x(F_i\setminus\{i\})).
\end{equation*}
Assuming that $j \in V^\F(i)$ and that configuration $x(F)$ maximizes ($\ref{eq:threshold}$), we have that 
$$|\mathbb{E}[W] | \geq 2\epsilon.$$
Therefore to bound (\ref{eq:simplebound}) for $j \in V^\F(i)$, it is enough to have an upperbound for
\begin{equation*}
 \P\left( \left| \hat{D}_n(x,F_i,i,j) -  |\mathbb{E}[W] |\right|  \geq \epsilon \right).
\end{equation*}

To do this, we observe that
\begin{equation*}
\left| \hat{D}_n(x,F_i,i,j) -  |\mathbb{E}[W] |\right| \notag 
\leq \left| \hat{D}_n(x,F_i,i,j) -  |W| \right| + \left| |W| -  |\mathbb{E}[W] |  \right|. \notag 
\end{equation*}

Then, using inequality (\ref{eq:enquadra}) we have
\begin{equation}
\left| \hat{D}_n(x,F_i,i,j) -  |\mathbb{E}[W] |\right| 
\leq |R| + \left| W -  \mathbb{E}[W]  \right|. \label{eq:lowerinv}
\end{equation}

Now, by (\ref{eq:lowerinv})
\begin{align}
 & \P\left( \left| \hat{D}_n(x,F_i,i,j) -  |\mathbb{E}[W] |\right|  \geq \epsilon \right)\\
& \leq \P\left(|R| \geq \frac{1}{2}\epsilon \right) + \P\left(\left| W -  \mathbb{E}[W]  \right| \geq \frac{1}{2}\epsilon \right).
\end{align}

Note that $\mathbb{E}[R] = 0$, thus by Bernstein inequality
\begin{equation} \label{eq:under1}
\P\left(|R| \geq \frac{1}{2}\epsilon \right) \leq 2\exp\left(-\frac{3n\epsilon^2}{4(6v+\epsilon)} \right),
\end{equation}
where $v$ is the same in (\ref{eq:variance2}).
By Bernstein inequality also we have
\begin{equation} \label{eq:under2}
 \P\left(\left| W -  \mathbb{E}[W]  \right| \geq \frac{1}{2}\epsilon \right) \leq 2\exp\left(-\frac{3n\epsilon^2}{4(6v_1+\epsilon)} \right),
\end{equation}
where $v_1$ is the same in (\ref{eq:variance1}).

Combining (\ref{eq:under1}) and (\ref{eq:under2}) we have for $j \in V^\F(i)$
\begin{equation*}
  \mathbb{P}\left(j \notin \hat{V}_n(i)\right) \leq 4\exp\left(-\frac{n\epsilon^2}{8v+\frac{4}{3}\epsilon}\right).
\end{equation*}

From this, it follows that
\begin{equation}
  \mathbb{P}\left( \mathcal{U}^\F_n(i) \right) \leq 4|V^\F(i)|\exp\left(-\frac{n\epsilon^2}{8v+\frac{4}{3}\epsilon} \right). \label{eq:under}
\end{equation}

Adding (\ref{eq:over}) and (\ref{eq:under}) we conclude the proof of Theorem (\ref{teo:bounded}).

\section*{Proof of Theorem \ref{teo:dob}}

  Let $z,z' \in \{-1,+1\}^{S}$ be two fixed configurations. For $i \in S$, $F_i \in \F$, let $(Y^z_t, Y^{z',F_i}_t)$ be a discrete time
  Markov chain taking values on $\{-1,+1\}^{F_i}$ with the following features.

  \begin{enumerate}
  \item The Ising model on $\{-1,+1\}^{F_i}$ with pairwise potential $J$ and boundary condition $z (F^c_i)$ is reversible with respect to the first marginal $Y^z_t$. 
  \item The Ising model on $\{-1,+1\}^{F_i}$ with pairwise potential $J^\F$ and  boundary condition $z'(F^c_i)$ is reversible with respect to the second marginal $Y^{z',\F}_t$. 
  \item The coupling chain $(Y^z_t, Y^{z',\F}_t)$ is irreducible and aperiodic, and has an unique invariant probability measure. Taking into the account items (1) and (2), this unique invariant probability measure is a coupling between the Ising models on 
  $\{-1,+1\}^{F_i}$ with interaction potentials $J$ and $J^\F$ and boundary conditions $z (F^c_i)$ and $z'(F^c_i)$ respectively. 
  \end{enumerate}
  
    We now construct $(Y^z_t, Y^{z',\F}_t)$ with $ t \in \mathbb{N}$. This can be done as follows. Let  $(I_t )_{t\ge 1}$ be an independent sequence of random variables uniformly distributed on $F_i$. For any $j \in F_i$ and $y\in \{-1,+1\}^{F_i}$, let also the probabilities $p_j(\cdot \;|\;y)$ and
$p^\F_j(\cdot \;| \;y)$ on $\{-1,+1\}$ be defined as follows. 
\[
p_j(+1\;|\;y)=\left \{1+e^{-2\sum_{k\in F_j}J(j,k)y(k)-2\sum_{k\notin
   F_i}J(j,k)z(k)}\right \}^{-1}\, ,
\]
\[
p^\F_j(+1\;|\;y)=\left \{1+e^{-2\sum_{k\in
    F_j}J^\F(j,k)y(k)-2\sum_{k\notin F_j}J^\F(j,k)z'(k)}\right \}^{-1}\,
.
\]

For any pair $(y, y')\in \{-1,+1\}^{F_i}\times \{-1,+1\}^{F_i}$, let $(\xi_t^{j, y, y'})_{t\geq 1}$, be an i.i.d. sequence of random variables taking values on  $\{-1,+1\}^2$ with distribution
\begin{align}
&\P\left(\xi_t^{j, y, y'}=(s,s)\right)=\min\left\{p_j(s|y), p^\F_j(s|y') \right\} , \notag\\
&\P\left( \xi_t^{j, y, y'}=(s,-s)\right)=\max\left\{p_j(s|y)-p^\F_j(s|y'), 0 \right\} ,\label{eq:couplenondiag}
\end{align}
for any $s \in \{-1,+1\}$.

Finally, let us assume that the sequences 
$(I_t )_{t\ge 1}$ and 
$(\xi_t^{j ,y, y'})_{t \ge 1}$, with $(y,y')\in \{-1,+1\}^{2F_i}$ and $j\in F_i$
 are all independent. The Markov chain $(Y^z_t, Y^{z',\F}_t)$ is constructed as follows.
For any $t\ge 1$ and any $j \in F$
\begin{equation*}
(Y^z_t(j), Y^{z',\F}_t(j))=(Y^z_{t-1}(j), Y^{z',\F}_{t-1}(j))\, ,\, \text{\;if\;} j\neq I_t\, 
\end{equation*}
and
\begin{equation}
(Y^z_t(j), Y^{z',\F}_t(j))=\xi_{t-1}^{j,Y^z_{t-1}, Y^{z',\F}_{t-1}}\, ,\, \text{\;if\;} j= I_t\;. \label{eq:couplingp}
\end{equation}

We stress the fact that the probabilities $p_j(\cdot|y)$ and
$p_j^\F(\cdot|y)$ depend on the fixed configurations $z(F^c_i) \in
\{-1,+1\}^{F^c_i}$ and $z'(F^c_i)\in \{-1,+1\}^{F^c_i}$ respectively. As a
consequence, the law of the Markov chain $(Y^z_t, Y^{z',\F}_t)$ depends on the
pair of fixed configurations $(z(F^c_i),z'(F^c_i))$.  Therefore,
a more explicit notation should mention all these details. This would
produce cumbersome things like $p_j(s|y(F_i\setminus {j}), z(F^{c}_i))$,
$p^\F_j(s|y(F_i\setminus {j}), z'(F^{c}_i))$. Hence we decided to use a simplified notation $p_j(\cdot|y)$ and
$p_j^\F(\cdot|y)$, respectively.

Let us assume that the initial value $(Y_0^z, Y_0^{z',\F})$ of the
chain is chosen according to its unique invariant probability measure.  For
every integer $t\geq 1$ we have
\begin{align}
  \P\left(Y^z_t(i) \neq Y^{z',\F}_t(i)\right)&= \P\left(Y^z_t(i) \neq Y^{z',\F}_t(i)\;,\;  I_t \neq i  \right )\notag \\
  &+\P\left(Y^z_t(i) \neq Y^{z',\F}_t(i)\;,\; I_t =
      i \right ).\label{eq:DobSeparate}
\end{align}

For the first term in the left hand side of the above equation we have
\begin{equation}
\P\left( Y^z_t(i) \neq Y^{z',\F}_t(i)\;,\; I_t \neq i \right )=\frac{|F_i|-1}{|F_i|}\P\left(Y^z_{t-1}(i) \neq Y^{z',\F}_{t-1}(i) \right). \label{eq:firstsepdob}
\end{equation}
 Substituting (\ref{eq:firstsepdob})  in (\ref{eq:DobSeparate}), and using the fact that the Markov chain is stationary, we obtain
\begin{equation}
\frac{1}{|F_i|}\P\left(Y^z_t(i) \neq Y^{z',\F}_t(i) \right)= \P\left(Y^z_t(i) \neq Y^{z',\F}_t(i)\;,\; I_t = i\right ).\label{eq:DobDivide}
\end{equation}

Now, we have
\begin{align}
&\P\left(Y^z_t(i) \neq Y^{z',\F}_t(i)\,,\,I_t = i \right)\notag \\
&=\P\left(Y^z_t(i) \neq Y^{z',\F}_t(i)\;,\;I_t = i \;,\;  Y^z_{t-1}(j) = Y^{z',\F}_{t-1}(j)\;\; \text{for all}\;\; j\in F_i \right ) \notag \\
&+ \P\left(Y^z_t(i) \neq Y^{z',\F}_t(i)\;,\;I_t = i  \; ,\; Y^z_{t-1}(j) \neq Y^{z',\F}_{t-1}(j) \;\text{for some} \;j\in F_i  \right ).\label{eq:dobsep}
\end{align}

Using  (\ref{eq:couplenondiag}) and (\ref{eq:couplingp}), the first term in the right hand side of  (\ref{eq:dobsep}) is bounded above by 
\begin{equation}
2\sup_{y\in \{-1,1\}^{F_i}}  \left (p_i(s|y)-p^\F_i(s|y) \right) \P\left( I_t = i\;,\; Y^z_{t-1}(j) = Y^{z',\F}_{t-1}(j)\;\; \text{for all}\;\; j\in F_i \right). \label{eq:dobestimate1}
\end{equation}

Using the mean value theorem, we have
\begin{equation}
\sup_{y\in \{-1,1\}^{F_i}} \left (p_i(s|y)-p^\F_i(s|y) \right) \leq \frac{1}{2}\sum_{l \notin F_i}|J(i,l)|.\label{eq:dobestimate2}
\end{equation}
Therefore an upperbound for expression (\ref{eq:dobestimate1}) is given by
\begin{align}
 &\frac{1}{|F_i|}\sum_{l \notin F_i}|J(i,l)| \P\left(Y^z_{t-1}(j) = Y^{z',\F}_{t-1}(j)\;\; \text{for all}\;\; j\in F_i\right)\notag\\
 &=\frac{1}{|F_i|}\sum_{l \notin F_i}|J(i,l)| \left[1 - \P\left(  Y^z_{t-1}(j) \neq Y^{z',\F}_{t-1}(j) \; \text{for some}\; j \in F_i \right )\right] \label{eq:outside2}
\end{align}

Let now study the second term of the right hand side of (\ref{eq:dobsep}). We first rewrite it as
\begin{equation*}
\sum_{\substack{U \subset F_i\\U\neq \emptyset }}\P\left(Y^z_t(i) \neq Y^{z',\F}_t(i)\;,\; I_t = i\;,\;\bigcap_{j \in U} \{Y^z_{t-1}(j) \neq Y^{z',\F}_{t-1}(j)\}\;,\; \bigcap_{j \in F_i \setminus U} \{Y^z_{t-1}(j) = Y^{z',\F}_{t-1}(j)\} \right ). \label{eq:111}
\end{equation*}

Therefore, proceeding as in (\ref{eq:dobestimate1}) and (\ref{eq:dobestimate2}), we obtain the following upperbound for the second term in the right hand side of (\ref{eq:dobsep}) 
\begin{align}
&\frac{1}{|F_i|} \sum_{\substack{U \subset F_i\\U\neq \emptyset }}\sum_{l \notin F_i}|J(i,l)|\P\left( \bigcap_{j \in U} \{Y^z_{t-1}(j) \neq Y^{z',\F}_{t-1}(j)\}\; , \;\bigcap_{k\in F_i\setminus U} \{Y^z_{t-1}(k) = Y^{z',\F}_{t-1}(k)\}\right)\notag \\
&+ \frac{1}{|F_i|}\sum_{\substack{U \subset F_i\\U\neq \emptyset }}\sum_{l \in U}|J(i,l)|\P\left(  \bigcap_{j \in U} \{Y^z_{t-1}(j) \neq Y^{z',\F}_{t-1}(j)\}\; , \; \bigcap_{k\in F_i\setminus U} \{Y^z_{t-1}(k) = Y^{z',\F}_{t-1}(k)\}\right). \label{eq:sum1}
\end{align}

The first part of  (\ref{eq:sum1}) can be rewritten as
\begin{equation}
\frac{1}{|F_i|}\sum_{l \notin F_i}|J(i,l)|\P\left(Y^z_{t-1}(j) \neq Y^{z',\F}_{t-1}(j) \; \text{for some}\; j \in F_i  \right ).\label{eq:firstdob}
\end{equation}

The second part of (\ref{eq:sum1}) can be rewritten as
\begin{equation*}
\frac{1}{|F_i|}\sum_{l \in F_i}\sum_{U \subset F_i: l \in U}|J(i,l)|\P\left( \bigcap_{j \in U} \{Y^z_{t-1}(j) \neq Y^{z',\F}_{t-1}(j)\}\bigcap_{k\in F_i\setminus U} \{Y^z_{t-1}(k) = Y^{z',\F}_{t-1}(k)\}\right)
\end{equation*}
and this is equal to
\begin{equation}
\frac{1}{|F_i|}\sum_{l \in F_i}|J(i,l)|\P\left( Y^z_{t-1}(l) \neq Y^{z',\F}_{t-1}(l) \right). \label{eq:seconddob2}
\end{equation}

Collecting together  (\ref{eq:dobsep}),  (\ref{eq:outside2}),  (\ref{eq:firstdob}), (\ref{eq:seconddob2}), we finally get the upperbound
\begin{equation}
\P\left(Y^z_t(i) \neq Y^{z',\F}_t(i) \right) \leq \sum_{l \notin F_i}|J(i,l)| +\sum_{l \in F_i}|J(i,l)|\P\left(Y^z_{t-1}(l) \neq Y^{z',\F}_{t-1}(l) \right). \label{eq:almostlastupper}
 \end{equation}
 
 To conclude the proof of the theorem, let  $Z$ and $Z'$ be two independent copies of the Ising models on 
 $\{-1,+1\}^{S}$ with potentials $J$ and $J^\F$, respectively. For a fixed realization of the pair $Z$ and $Z'$, construct as before the coupled chains $(Y^Z_t, Y^{Z',\F}_t)$ taking values on $\{-1,+1\}^{2F_i}$, and having  $Z(F^c_i)$ and $Z'(F^c_i)$ as boundary conditions.
 
Using inequality (\ref{eq:almostlastupper}) and taking the expectation with respect to $(Z,Z')$, we have
\begin{equation} \label{eq:inefinal}
\mathbb{E}\left[\P\left(Y^Z_t(i) \neq Y^{Z', \F}_t(i)  \right)\right] \leq \sum_{l \notin F_i}|J(i,l)| +\sum_{l \in F_i}|J(i,l)|\mathbb{E}\left[\P\left(Y^Z_{t-1}(l) \neq Y^{Z', \F}_{t-1}(l) \right)\right].
 \end{equation}
  Now observe that 
\begin{equation*}
\mathbb{E}\left[\P\left(Y^Z_t(j) \neq Y^{Z', \F}_t(j)  \right)\right]  = \P\left(Y(j) \neq Y^\F(j) \right),
 \end{equation*}
 for any $j \in F$, where  $Y(j)$ and $Y^\F(j)$ are the projections on site $j$ of realizations of the Ising model with pairwise potential $J$ and $J^\F$, respectively.    
 From this identity and inequality (\ref{eq:inefinal}), it follows that
 \begin{equation*}
 \P\left(Y(i) \neq Y^\F(i) \right) \leq \sum_{l \notin F_i}|J(i,l)| +\sum_{l \in F_i}|J(i,l)|\P\left(Y(l) \neq Y^\F(l)\right).
 \end{equation*}
Finally, taking the supremum for all $i \in F$ we have
\begin{equation*}
\sup_{i \in S}\P\left(Y(i) \neq Y^\F(i) \right) \leq \sup_{i\in S}\sum_{l \notin F_i}|J(i,l)| +r\sup _{i\in S}\P\left( Y(i) \neq Y^\F(i)\right),
 \end{equation*} 
  which concludes the proof.
 
%
%\newpage
%\noindent {\bf \Large Appendix 2}
%\section*{Pseudocode of the algorithm Sample Generator}
%
%
%
%\begin{algorithm}[h]  \label{alg:sampler}
%{\bf \normalsize Algorithm}{\;\;Sample generator}
%\begin{algorithmic}[1]
%\FOR {$k: 1 \to n$}
%\STATE {$t=0$ and $Y_t(i) = +1$ and $Y'_t(i)=-1$ for all $i \in B^L(0)$.}
%\WHILE {$Y_t \neq Y'_t$} 
%\STATE $t \gets t+1$
%\STATE {Generate a random variable $I_t $ uniformly distributed on $B^L(0)$.}
%\FOR {$j \in B^L(0)$}
%\IF{ $j=I_t$ }
%\STATE {Choose $(s,s') \in \{-1,1\}^2$ with probability 
%\begin{align*}
%&P_j\left(s,s' | y, y' \right)=\min\left\{p_j(s|y), p_j(s|y') \right\}, \; \text{if} \; s=s' \notag\\
%&P_j\left(s,s' | y, y' \right)=\max\left\{p_j(s|y)-p_j(s|y'), 0 \right\}, \; \text{if} \; s\neq s'  
%\end{align*}
%where for any $y\in \{-1,+1\}^{B^L(0)}$ is given by 
%
%\[
%p_j(s\;|\;y)=\left \{1+e^{-2s\sum_{k\in B^L(0)}J(j,k)y(k)}\right \}^{-1}.
%\]
%}
%\STATE  $Y_t(j) \gets s$ and $Y'_t(j) \gets s'$, 
%\ELSE
%\STATE $Y_t(j) \gets Y_{t-1}(j)$ and $Y'_t(j) \gets Y'_{t-1}(j)$
%\ENDIF
%\ENDFOR
%\ENDWHILE
%\STATE $X_k=Y_t$
%\ENDFOR
%\RETURN Sample $X_1, \ldots, X_n$.
%\end{algorithmic}
%\end{algorithm}

\vspace{0.5cm}
\noindent {\bf Acknowledgements}

This article was produced as part of the activities of FAPESP Center for Neuromathematics (FAPESP grant 2013/ 07699-0).  This work is part of USP project ``Mathematics, computation, language, and the brain", and CNPq project ``Stochastic modeling of the brain activity" (grant 480108/2012-9). AG was partially supported by a CNPq fellowship
(grant 309501/2011-3). DYT was supported by FAPESP grant 2008/08171-0 and Pew Latin American Fellowship. EO was partially supported by Prin07:20078XYHYS.  DYT thanks the hospitality of Universit\`a di Roma Tre, Dipartimento di Matematica. EO thanks the hospitality of NUMEC, USP.

The authors would like to thank an anonymous referee who pointed an error in the statement of Theorem 3 in a earlier version of the article. The authors also thank Luiz Lana and Matthieu Lerasle for many illuminating discussions.

\bibliographystyle{plainnat}
\bibliography{Identify_Ising}

\vskip30pt

Antonio Galves

Instituto de Matem\'atica e Estat\'{\i}stica

Universidade de S\~ao Paulo

Caixa Postal 66281

05315-970 S\~ao Paulo, Brasil

e-mail: {\tt galves@usp.br}
\bigskip

  Enza Orlandi
 
Dipartimento di Matematica

Universit\`a  di Roma Tre

 L.go S.Murialdo 1, 00146 Roma,  Italy. 

email: {\tt orlandi@mat.uniroma3.it}
\bigskip

Daniel Yasumasa Takahashi

Institute of Neuroscience

Princeton University

Princeton, 08648, USA

e-mail: {\tt takahashiyd@gmail.com}

\end{document}